\documentclass[12pt]{amsart}
\usepackage[margin=1.0in]{geometry}
\usepackage{mathptmx}
\usepackage[foot]{amsaddr}
\begin{document}
\title[The Banach-Tarski Paradox in Bely's Modernist Novel]{From Poland to \emph{Petersburg}:\\
The Banach-Tarski Paradox in Bely's Modernist Novel}

\author{Noah Giansiracusa$^*$} 
\address{$^*$Assistant Professor of Mathematics, Swarthmore College} 
\author{Anastasia Vasilyeva$^{**}$} 
\address{$^{**}$Undergraduate Economics Major, Swarthmore College}
\email{ngiansi1@swarthmore.edu, avasily1@swarthmore.edu}

\begin{abstract}
Andrei Bely's novel \emph{Petersburg}, first published in 1913, was declared by Vladimir Nabokov one of the four greatest masterpieces of 20th-century prose.  The Banach-Tarski Paradox, published in 1924, is one of the most striking and well-known results in 20th-century mathematics.  In this paper we explore a potential connection between these two landmark works, based on various interactions with the Moscow Mathematical School and passages in the novel itself.
\end{abstract} 

\maketitle

\section{Introduction}
Andrei Bely (1880-1934) was a poet, novelist, and theoretician who helped lead the Symbolist movement in Russia in the early 20th century.  His most regarded work is the modernist novel \emph{Petersburg}, published serially in 1913-1914 then in a revised and shortened form in 1922.  Vladimir Nabokov (1899-1977) famously declared \emph{Petersburg} one of the four greatest works of prose of the 20th century, along with Joyce's \emph{Ulysses}, Kafka's \emph{Metamorphosis}, and Proust's \emph{\`A la recherche du temps perdu}.  As with these three other works, the scholarly literature on \emph{Petersburg} is vast and continues to grow; for instance, a volume of essays on \emph{Petersburg} celebrating its centennial was just released \cite{Centennial}.  

An emerging direction of scholarship has been to better understand the origins and meaning of the remarkably frequent mathematical imagery in \emph{Petersburg} \cite{Szilard1,Szilard2,Naming,Svetlikova,Columbia,Noah1}.  That Bely had an interest in math is not surprising: his father was the influential mathematician Nikolai Bugaev\footnote{It was to protect the reputation of his famous father that Boris Nikolaevich Bugaev, better known as Andrei Bely, chose to publish under a pseudonym \cite[p. 31]{Mochulsky}.} (1837-1903) who is credited with creating the Moscow Mathematical School, one of the most active and successful groups of mathematicians in recent history \cite{Demidov,Demidov-150}.  Moreover, while Bely was studying the natural sciences at Moscow University, prior to focusing professionally on the literary arts, his close classmates included two students of Bugaev: Nikolai Luzin (1883-1950), who took the reins of the Moscow Mathematical School for many years, and Pavel Florensky (1882-1937), a fellow Symbolist who was a Russian Orthodox priest and a mathematician (among other things \cite{Florensky}).  The scholarship cited above looks at these mathematicians and others in their circle and explores how their philosophical ideas, mystical inclinations, and political affiliations help us better understand the complex, multi-layered Symbolist novel \emph{Petersburg}.

One of the most ubiquitous and significant symbols in \emph{Petersburg} is an expanding, and ultimately exploding, sphere.  This describes the sensation three characters feel in their chest from symptoms of anxiety and heart disease, and it foreshadows a bomb that one of them gives to another to use on the third.  Various scholars have argued that the sphere plays the role in the novel of a Wagnerian leitmotif and that its metaphysical import is cosmic thermodynamics, a Dionysian leap, circularity of thought/evolution/behavior, etc.  As we explain in the present paper, there is also a previously unnoticed \emph{mathematical} interpretation of Bely's expanding spheres, but it is based on a theorem, the astonishing Banach-Tarski Paradox, that appeared a decade after the novel.  

Could Bely have inspired Stefan Banach (1892-1945) and Alfred Tarski (1901-1983)?  Could earlier versions of their theorem have inspired Bely?  While the answer to both questions is ``probably not,'' we nonetheless trace the historical developments that could potentially have put Bely in contact with the mathematical ideas leading up to Banach-Tarski.  Bely believed in spiritual connections and mystical predictions, so in this vein we also explore the (sometimes startling) coincidences uniting \emph{Petersburg} with the Banach-Tarski Paradox.  

This paper is thus the story, part history and part mystery, of an unlikely link between math and literature.  We start by recalling the Banach-Tarski Paradox from a mathematical perspective.  Next, we trace the history of this theorem and discuss the mathematical landscape that led to it.  Banach and Tarski were part of a Polish School of Mathematics, and we present a bridge between Andrei Bely and this school stemming from the geographic displacement of Waclaw Sierpi\'nski (1882-1969) brought upon by World War I.  Having set the stage historically, we then turn to the novel \emph{Petersburg} itself by analyzing and contextualizing several passages in terms of the Banach-Tarski Paradox.  We conclude by discussing some predictions and coincidences surrounding Bely that other scholars have noted and we place the Banach-Tarski Paradox in this context.

\section{The Banach-Tarski Paradox}

The strong form of the Banach-Tarski Paradox, which appeared in their original paper \cite{Banach-Tarski}, is the following:
Given any two bounded sets $A,B\subset \mathbb{R}^n$ with nonempty interior and $n\ge 3$, it is possible to decompose $A=\sqcup_{i=1}^k A_i$ into a finite number of disjoint subsets $A_i$, and similarly for $B=\sqcup_{i=1}^k B_i$, such that for each $i=1,\ldots,k$ there exists a Euclidean isometry sending $A_i$ onto $B_i$.  The two most commonly stated special cases of this theorem are when $A\subset \mathbb{R}^3$ is a solid ball and $B$ is either two identical, disjoint copies of $A$ (in which case $k=5$ suffices) or a solid ball of any radius larger than that of $A$.  The paradoxical aspect, of course, is that intuitively it should not be possible to increase the volume of an object simply by cutting it into a few pieces then translating, rotating, and possibly flipping those pieces---but some of the subsets $A_i,B_i$ are necessarily not Lebesgue measurable, so volume is a tricky concept here and the ``cuts'' one needs are not possible in reality.  

To get a vague sense of how the Banach-Tarski Paradox works, it is helpful to consider the following much simplified variant which is a circular version of the famous infinite-room hotel of David Hilbert (1862-1943).  Fix a point $x_0\in S^1$ on the unit circle and consider the infinite sequence of points $X = \{x_n\}_{n=0}^\infty \subset S^1$ where $x_n$ is the result of rotating $x_0$ counterclockwise by $n$ radians.  All the points of this sequence are distinct since $2\pi$ is an irrational number.  For any integer $m\ge 1$ we can rotate $X$ clockwise by $m$ radians to obtain the set \[X\sqcup \{x_{-1}\} \sqcup \cdots \sqcup \{x_{-m}\}.\]  Thus, by rotating this set $X$ we are able to add any finite number of points to it.  Increasing the Lebesgue volume of a set requires adding an uncountable number of points, so we are still a long way from the Banach-Tarski Paradox, but at least we see here how isometries don't necessarily behave how one would expect when infinite cardinalities are involved.  

One of the main steps of the actual Banach-Tarski Paradox is to show that the special orthogonal group $SO(3)$ contains a subgroup isomorphic to the free group on two generators.  However, for the purposes of the present paper we do not need to understand any of these technical details.  Since there are many excellent expository accounts of the Banach-Tarski Paradox and its proof, we therefore simply refer the interested reader to these \cite{PeaSun,Wagon,Weston,Raman,Tao}.

\section{Tracing the history}

The discovery of non-Euclidean geometries in the 19th century was a tremendous surprise and a paradigm leap in mathematics; but the Banach-Tarski Paradox upends our most basic understanding of Euclidean space itself.  In this section we look at where it came from, who was involved, and how Andrei Bely fits in.

\subsection{Setting the stage for Banach-Tarski}

The stunning paper of Banach and Tarski \cite{Banach-Tarski}, which came out two years after Bely's revision of \emph{Petersburg}, builds on three important earlier enigmatic concepts/constructions: the Axiom of Choice, introduced by Ernst Zermelo (1871-1953) in 1904 \cite{Zermelo}; the first example of a subset of $\mathbb{R}$ that is not Lebesgue measurable, introduced by Giuseppe Vitali (1875-1932) in 1905 \cite{Vitali}; and the Hausdorff Paradox from 1914 \cite{Hausdorff}, the year \emph{Petersburg} completed its publication in serial format.  This paradox of Felix Hausdorff (1868-1942) informally stated is that, if we are willing to ignore a countable set of points, then there is a Euclidean isometry mapping half of a sphere in $\mathbb{R}^3$ onto a third of the same sphere.  

Another important piece of the mathematical history leading up to the Banach-Tarski Paradox is a paper, also from 1914, in which Stefan Mazurkiewicz (1888-1945) and Sierpi\'nski construct a subset of the plane $\mathbb{R}^2$ that is carried by a Euclidean isometry onto half of itself, roughly speaking \cite{Mazurkiewicz-Sierpinski}.  Sierpi\'nski continued working on mathematics related to these topics for many years; in 1948 he published an alternative proof of the Banach-Tarski Theorem based on the Hausdorff Paradox \cite{Sierpinski}.  

What led to this surge of interest in mathematical paradoxes?  With the introduction of set theory by Georg Cantor (1845-1918) in the late 19th century, the foundations of mathematics in the early 20th century were quickly being rewritten.  Functional analysis and measure theory particularly flourished in this period with their newfound rooting in set theory.  But this systematic recasting of mathematics brought deep logical issues ineluctably to the forefront, especially in light of the famous set-theoretical paradox of Bertrand Russell (1872-1970) published in 1901.  Mathematics could not move forward without a better understanding of set theory and its myriad consequences---and one way of navigating these uncharted waters was to uncover the most perplexing, counterintuitive results to test the limits of the subject.  

Unease from the abundance of mathematical strangeness that was emerging at the time led some to question whether the new foundations of mathematics had perhaps strayed beyond mathematics itself.  Since geometry is arguably the field of mathematics most closely built upon our intuition, a particularly fertile ground for exploring set-theoretical paradoxes was in the geometric aspects of measure theory.  Thus the Banach-Tarski Paradox and the geometric paradoxes leading up to it were in essence a clarion call to mathematicians world-wide: if you accept the new set-theoretic foundations of mathematics, then you must also accept a geometry that no longer accedes to your intuition.  

\subsection{Mathematical schools in the early 20th century}

To begin tracing the possible interactions between Bely and the Banach-Tarski Paradox, we need to step back and see how different geographical schools of mathematics were forming, or solidifying, in the early 20th century in response to the Cantorian set-theoretic revolution.

The book \cite{Naming} tells a fascinating story of how French mathematicians in the late 19th and early 20th centuries accepted the mathematical content of Cantor's work but, true to their avowed secularism, avoided the more philosophical and religious interpretations that were arising.  In stark contrast, the Moscow Mathematical School led by Bely's father, Bugaev, and his students viewed math as underlying the universe and all realms of human knowledge (notably psychology, history, and aesthetics, not just the more typical natural sciences), and in doing so enthusiastically embraced what they saw as the mystical dimensions of set theory.  

An amusing exception to this dichotomy is the following.  In 1905 the French mathematicians \'Emile Borel (1871-1956), Ren\'e-Louis Baire (1874-1932), Henri Lebesgue (1875-1941), and Jacques Hadamard (1865-1963) were discussing Zermelo's recently introduced Axiom of Choice.  The first three of them all rejected the axiom, whereas Hadamard said \cite[57]{Naming}: ``the question of what is a correspondence that can be described is a matter of psychology and relates to a property of the mind outside the domain of mathematics.'' To consider a mathematical issue under the purview of psychology would fit right in with the Moscow school, though to consider psychology outside the domain of mathematics certainly would not!

Another important group of mathematicians making amazing progress exploring the new set-theoretic foundations was the Polish School of Mathematics (consisting of the Lw\'ow, Warsaw, and Krak\'ow chapters).  Among many others, this school included Banach and Tarski, as well as two of the mathematicians mentioned above in the lead-up to the Banach-Tarski Paradox: Mazurkiewicz and, crucially for this story (as we soon discuss), Sierpi\'nski.  Throughout the 20th century, and even to this day, many of the results on paradoxical geometric constructions come from Polish mathematicians.  One can consult Wagon's book \cite{Wagon}, particularly the historical notes and extensive bibliography, to get a sense of this national mathematical cohesion.  

Thus, at the beginning of the 20th century we have the French atheists working on measure theory and analysis; the Moscow mystics in response working on these same subjects while involving them in grandiose theories of the universe; and the Polish paradox creators probing the counterintuitive depths of these same subjects.  

\subsection{Connecting the dots}

To connect the author of \emph{Petersburg}, Andrei Bely, to the Moscow Mathematical School is easy: as we mentioned earlier, Bely's father Bugaev founded the school and while at university Bely was friends with Luzin and Florensky, two math students of Bugaev who played important roles in this school.  Evidently Luzin and Florensky debated politics, philosophy, and mathematics with Bely \cite[p. 79]{Naming}.  Since these three graduated between 1903 and 1905, it is safe to assume that Bely was aware of the Cantorian tumult dominating mathematics at that time and of the Moscow school's response to it.

But then to connect Bely to the Banach-Tarski Paradox and its precursors we need a link between the Moscow school and the Polish school.  In fact, such a link did exist and was provided by Sierpi\'nski.  In 1914, while he was a professor at the university of Lw\'ow in Poland (and the year his planar paradox with Mazurkiewicz was published), Sierpi\'nski went on a trip to Russia to visit some family members.  World War I broke out at this time and he was imprisoned by the Tsarist authorities in Vyatka (now Kirov), Russia.  Luzin was aware of  Sierpi\'nski's mathematical prowess, and upon receiving news of his imprisonment Luzin managed to have him released and transferred to Moscow \cite{Art}.   Sierpi\'nski spent the remaining war years 1915-1918 deeply immersed in the activities of the Moscow Mathematical School, working closely with Luzin, before returning to Poland to spend the rest of his career at the University of Warsaw \cite{Demidov}.  

Luzin, therefore, is a direct connection between his classmate and friend, the author Andrei Bely, and his war-time collaborator, the mathematician Sierpi\'nski who would go on to become an expert in Banach-Tarski and who already by his Moscow years was an established expert on geometric paradoxes.  But the problem, of course, is that Sierpi\'nski arrived in Moscow the same year that Bely finished publishing his novel \emph{Petersburg}.  If the expanding/bursting sphere symbol had been an addition to the 1922 edition, absent from the 1913-1914 serial publication, then one could somewhat confidently conjecture that Sierpi\'nski, through Luzin, had influenced Bely's writing---but this is \emph{not} the case: Bely's Banach-Tarski-esque spheres are just as prominent and ubiquitous in the first version of his novel as they are in the revision.  

There is just a small window left open here.  The fact that Luzin went to such efforts to rescue and relocate Sierpi\'nski in 1914 suggests that Luzin was at that time already somewhat familiar with Sierpi\'nski's work, which means he may have been following Sierpi\'nski for several years---so Bely may have been as well, or at least he may have gotten wind of the types of results Sierpi\'nski and others in the Polish school were proving while he was writing his novel.  But the geometric paradoxes predating 1914 just don't bear much of a resemblance to the spherical symbolism in \emph{Petersburg}.  Thus the probability of contemporary geometric paradoxes influencing Bely's novel, while distinctly nonzero, seems fairly low.  And conversely, Bely likely was aware of earlier set-theoretical issues such as Russell's Paradox and the Axiom of Choice, but the momentum toward uncovering geometric paradoxes in the Polish school was already strongly established and the Banach-Tarski Paradox builds directly off of earlier mathematical results, so it is unlikely Bely's literary work influenced these mathematicians.  

And so, without further evidence revealing an explicit influence in one direction or the other, we must accept the likelihood that the similarities (discussed in the next section) between the expanding spheres in the work of Bely and Banach-Tarski are merely a coincidence, albeit a rather magnificent one---even more so, as we discuss in the final section of this paper, due to a happenstance concerning the colloquial name by which the Banach-Tarski Paradox eventually came to be referred.  

\section{Expanding Spheres in \emph{Petersburg}}

The Banach-Tarski Paradox, in one of its common formulations, involves a sphere that increases in size, and it accomplishes this mysterious task by breaking apart into pieces that are rotated into place.  There are many passages in \emph{Petersburg}\footnote{Note: All the \emph{Petersburg} quotes in the rest of this paper are from the Macguire-Malmstad English translation of the 1922 version \cite{Petersburg}, but in most cases similar passages exist in the 1913-1914 edition as well.} that seem reminiscent of this situation.  For instance \cite[p. 14]{Petersburg}: ``His heart pounded and expanded, while in his breast arose the sensation of a crimson sphere about to burst into pieces.''  And \cite[p. 284]{Petersburg}: ``his soul was becoming the surface of a huge, rapidly growing bubble, which had swollen into Saturn's orbit.  Oh, oh, oh! Chills ran through Nikolai Apollonovich.  Winds wafted into his forehead.  Everything---was bursting.''  Next, an even more elaborate, and nightmarish, example:
\begin{quote}
A bomb is a rapid expansion of gases.  The sphericality of the expansion evoked in him a primordial terror, long forgotten.  In his childhood he had been subject to delirium.  In the night, a little elastic blob would sometimes materialize before him and bounce about---made perhaps of rubber, perhaps of the matter of very strange worlds.  [...]  Bloating horribly, it would often assume the form of a spherical fat fellow.  This fat fellow, having become a harassing sphere, kept on expanding, expanding, and expanding and threatened to come crashing down upon him. [...]  And it would burst into pieces.  Nikolenka would start shrieking nonsensical things: that he too was becoming spherical, that he was a zero, that everything in him was zeroing---zeroing---zero---o---o... \cite[pp. 157--158]{Petersburg}
\end{quote}
Nikolenka was hiding a bomb in his room, waiting to use it on his father as part of a terroristic revolutionary plot that he had reluctantly involved himself in.  The expanding sphere represents this unstoppable bomb and also Nikolenka's internal anxiety over his impending patricidal act.  

The Banach-Tarski Paradox (which, as we have discussed, did not exist at the time of Bely's writing) reached many in the mathematical community as a sort of anxiety-inducing bomb as well: if such an illogical conclusion, defying our most innate intuitions about shape and volume, follows from the basic axioms underlying all of mathematics, then something might be wrong with these axioms and hence with the entire edifice of modern mathematics.  For those in the early 20th century who had dedicated their careers to solidifying the set-theoretic foundations of mathematics, the Banach-Tarski Paradox could well have caused nightmarish visions like the one Bely describes.   

Yuri Manin (born 1937), who was a big part the Moscow Mathematical School in its later years, mentioned these issues in an interview not too long ago:
\begin{quote}
Take for instance the theorem of Banach-Tarski. You start with a ball, and it turns out that you can cut it into five pieces, rearrange them, put them back together, and you obtain two balls of the same size as the initial one.  This construction tells us a lot. For example, to the critics of the set-theoretic approach in general, it means that if this view leads one to such an assertion, then it is not mathematics, but some sort of wild nonsense. For logicians it is an example of a paradoxical application of the axiom of choice of Zermelo and so an argument against accepting it. And aside from all this, it is very beautiful geometry. \cite{Manin}
\end{quote}

Manin goes on to explain how Banach-Tarski, and other apparent paradoxes, did not cause the foundations of mathematics to crumble, nor even decelerate mathematical progress, but instead helped us understand various phenomena on a deeper level:
\begin{quote}
Several such paradoxes were discovered during the time of transition between classical mathematics and set theoretic mathematics. There was the theorem that a curve could fill the square. There were many such things, and they taught us a lot.  Many people thought that this was pure fantasy, but the newly trained imagination allowed one to recognize ``paradoxical'' behavior of Fourier series, to understand Brownian motion, then to invent wavelets, and it turned out that these were not at all fantasies but almost applied mathematics. \cite{Manin}
\end{quote}
We must remember that \emph{Petersburg} was written during this time of transition that Manin mentions, so if it indeed encodes mathematical developments in its literary symbolism then we should expect the anxious, unsettled, and emotional initial reactions Manin describes rather than the calmer, enlightened, and rational hindsight that developed later in the 20th century.  

Manin also explains a useful way of understanding the Banach-Tarski construction intuitively, or at least of seeing why it may not be as baffling as it first appears:
\begin{quote}
The key point was that we must not imagine ``pieces'' as solid material objects, but rather clouds of points. We must imagine that a ball consists of indivisible points. You are allowed to call a ``piece'' any subset of these points, you can move it and turn it around, but only as a whole, moving it as a single object, so that the pairwise distances between points remain the same. So you split the sphere not into solid pieces, but into five clouds. And these clouds can mutually penetrate each other; in fact, there's nothing solid about them. They have no volume, no weight, they are wonderful objects of a highly trained imagination. [...] The message here is that if you make a dust of individual points out of your initial ball, there will be enough points to fill two, or three, or even an infinity of balls of arbitrary sizes. \cite{Manin}
\end{quote}
Curiously, another important symbol in \emph{Petersburg} is swarms of small indistinguishable objects, essentially clouds of points.  So Manin's explanation of Banach-Tarski seems, unintentionally and coincidentally, to unite two of Bely's main symbols in the novel. 

The following passage in \emph{Petersburg} conveys Manin's view of Banach-Tarski---swarms of points rotating to form any possible shape (the ``strong form'' of the paradox)---to an almost uncanny extent, especially considering that Manin's words come almost a hundred years after Bely's:
\begin{quote}
Ableukhov's eyes saw bright patches and dots of light, and iridescent dancing spots with spinning centers.  [...] And the misty spots and stars, like foam on bubbling blackness, would unexpectedly and suddenly form into a distinct picture: of a cross, a polyhedron, a swan, a light-filled pyramid.  And all would fly apart.  \cite[p. 93]{Petersburg}
\end{quote}

It is tempting to view these passages in \emph{Petersburg}---expanding spheres bursting into pieces, dancing spots with spinning centers forming various geometric shapes---as a deliberate literary manifestation of the Banach-Tarski Paradox,  But alas, even with the connection between Bely and the Polish school provided by Luzin and Sierpi\'nski, the chronology, as discussed earlier, simply does not support this interpretation.

This seems an appropriate time to mention the following quote from a biographer of Bely \cite[p. 51]{Elsworth}: ``He [Bely] came to the conviction that scientific discoveries tend to be preceded by artistic ones, that, for instance, it was the concrete realization of form in Greek sculpture that made possible the theorems of geometry.''  Whether causal or coincidental, Bely corroborates his own theory here by imbuing \emph{Petersburg} with symbolism suggestive of mathematics that would appear a decade later.  Strangely enough, there is a further connection between \emph{Petersburg} and the Banach-Tarski Paradox that, while quite astounding, is even more assuredly coincidental than the connections discussed in this section.

\section{Cosmic Coincidences}

Let us begin this final section by illustrating how eerily, and inexplicably, prescient Bely was at times.  One example of this phenomenon concerns \emph{Petersburg} itself, as one Bely critic explains:
\begin{quote}
Bely acknowledged in his posthumously published memoir that he had modeled Lippanchenko on the notorious double agent Evno Fishelevich Azef, who had worked concurrently for the Tsarist secret police and the terrorist arm of the Socialist Revolutionary Party. (Curiously, Bely gave his agent-provocateur a surname that resembled to an uncanny degree the pseudonym Azef used while living in Berlin, "Lipchenko"---a coincidence that subsequently amazed Bely because he could not have known the pseudonym when he was working on the novel.)  \cite[p. 144]{Alexandrov}
\end{quote}
Another example concerns Bely's death, which the literary scholar Mochulsky describes thusly: 
\begin{quote}
On July 17, 1933, in Koktebel, Bely suffered a sun stroke; he died in Moscow on January 8, 1934.  He had predicted his death in a poem of 1907: 
\begin{center}
I trusted the gold brilliance,\\
But died from the sun's arrows.\\
I measured the centuries in thought,\\
But could not live my own life well.
\end{center}
\cite[p. 226]{Mochulsky}
\end{quote}

The Sun is an important symbol in the occultist anthroposophy movement in which Bely was involved while writing \emph{Petersburg}.  Specifically, one of the names for Christ in anthroposophy is ``Sun Being'' \cite[p. 13]{Alexandrov}.  Bely once enthusiastically exclaimed \cite[p. 62]{Mochulsky}: ``She is already in our midst, with us, embodied, vital, near --- this finally recognized Muse of Russian Poetry has proven to be the Sun.''  And here in \emph{Petersburg} itself is a beautiful passage involving the Sun in a context that at least vaguely reminds one of Banach-Tarski:
\begin{quote}
There was no Earth, no Venus, no Mars, merely three revolving rings.  A fourth one had just blown up, and an enormous Sun was still preparing to become a world.  Nebulae whirled past.  Nikolai Apollonovich had been cast into measureless immensity, and distances flowed. \cite[p. 167]{Petersburg}
\end{quote}
Indeed, we have one spherical planet that has blown up and others that have been transformed into a different geometric shape; the nebulae are whirling clouds of dust and so quite reminiscent of Manin's description of Banach-Tarski; we even see mention of the words ``distances'' and ``measureless,'' which are two ingredients in Banach-Tarski if one recasts measureless as non-measurable.    But this cosmic coincidence truly takes form when one realizes that one of the popular names for the Banach-Tarski Paradox (which first appeared well after Banach-Tarksi's paper, and hence also \emph{Petersburg}) is the ``Pea and the Sun Paradox,'' since a solid sphere the size of a pea could be cut into pieces and reassembled into one the size of the Sun \cite{PeaSun}.

So, in summary, the expanding spheres throughout \emph{Petersburg}, bursting into pieces, into swirling point clouds, the ``measureless immensity'' mentioned at least six times in the book, we have interpreted all this in light of the Banach-Tarski Paradox.  And conversely, the Banach-Tarski Paradox, also known as the Pea and the Sun Paradox, in Bely's world takes on the mystical meaning of a tiny, insignificant object growing into Christ himself, or into the muse of Russian poetry.  As we have discussed, a causal relationship in either direction between \emph{Petersburg} and Banach-Tarski is unlikely, so we must simply accept these thematic and linguistic connections as cosmic coincidences.  Bely would be proud of himself for having achieved yet another inexplicable premonition---and for demonstrating that indeed artistic discoveries often predate, and impel, scientific ones.

Let us conclude with one final coincidence, a numerological one.  The number five often appears, sometimes conspicuously, in \emph{Petersburg}.  For instance:
\begin{itemize}
\item ``There the houses merged cubelike into a regular, five-story row.'' [p. 10]
\item ``Thus he had been sleeping for five years.'' [p. 15]
\item ``Then with all five fingers he began strumming the guitar.'' [p. 67]
\item ``He had been turning the wheel for a mere five years.'' [p. 232]
\end{itemize}
And most strikingly we see the number five oddly appearing in a passage that, once again, involves expansion and immeasurability:
\begin{quote}
I was growing, you see, into an immeasurable expanse, all objects were growing along with me, the room and the spire of Peter and Paul.  There was simply no place left to grow.  And at the end, at the termination---there seemed to be another beginning there, which was most preposterous and weird, perhaps because I lack an organ to grasp its meaning.  In place of the sense organs there was a ``zero.''  I was aware of something that wasn't even a zero, but a zero minus something, say five, for example. \cite[p. 182]{Petersburg}
\end{quote}
The prevalence of the number five has been interpreted in terms of anthroposophic numerology \cite{Columbia}, but we shall close this paper simply by recalling that five is also the number of pieces into which a sphere must be divided in order to double it via Banach-Tarski.

\bibliographystyle{alpha}
\bibliography{bib}

\end{document}